
\input amssym.def
\input amssym

\documentstyle{amsppt}
\magnification=1200

\font\cal=cmsy10
\font\small=cmr5

\hoffset 0,6truecm
\voffset 1truecm
\hsize 15,2truecm
\vsize 22truecm

\catcode`\@=11
\long\def\thanks#1\endthanks{%
  \ifx\thethanks@\empty \gdef\thethanks@{%
     \tenpoint#1} 
     \else \expandafter\gdef\expandafter\thethanks@\expandafter{%
     \thethanks@\endgraf#1}%
  \fi}
\def\keywords{\let\savedef@\keywords
  \def\keywords##1\endkeywords{\let\keywords\savedef@
  \toks@{\def\usualspace{{\it\enspace}}\tenpoint} 
  \toks@@{##1\unskip.}%
  \edef\thekeywords@{\the\toks@\frills@{%
  {\noexpand\it 
  Key words and phrases.\noexpand\enspace}}\the\toks@@}}%
  \nofrillscheck\keywords}
\def\makefootnote@#1#2{\insert\footins
 {\interlinepenalty\interfootnotelinepenalty
 \tenpoint\splittopskip\ht\strutbox\splitmaxdepth\dp\strutbox 
 \floatingpenalty\@MM\leftskip\z@\rightskip\z@
 \spaceskip\z@\xspaceskip\z@
 \leavevmode{#1}\footstrut\ignorespaces#2\unskip\lower\dp\strutbox
 \vbox to\dp\strutbox{}}}
\def\subjclass{\let\savedef@\subjclass
  \def\subjclass##1\endsubjclass{\let\subjclass\savedef@
  \toks@{\def\usualspace{{\rm\enspace}}\tenpoint}
  \toks@@{##1\unskip.}%
  \edef\thesubjclass@{\the\toks@\frills@{{%
  \noexpand\rm1991 {\noexpand\it Mathematics
  Subject Classification}.\noexpand\enspace}}%
  \the\toks@@}}%
\nofrillscheck\subjclass}
\def\abstract{\let\savedef@\abstract
 \def\abstract{\let\abstract\savedef@
  \setbox\abstractbox@\vbox\bgroup\noindent$$\vbox\bgroup
  \def\envir@end{\endabstract}\advance\hsize-2\indenti
  \def\usualspace{\enspace}\tenpoint \noindent 
  \frills@{{\smc Abstract.\enspace}}}%
 \nofrillscheck\abstract}

\outer\def\endtopmatter{\add@missing\endabstract
 \edef\next{\the\leftheadtoks}\ifx\next\empty
  \expandafter\leftheadtext\expandafter{\the\rightheadtoks}\fi
 \ifmonograph@\else
   \ifx\thesubjclass@\empty\else \makefootnote@{}%
        {\thesubjclass@}\fi
   \ifx\thekeywords@\empty\else \makefootnote@{}%
        {\thekeywords@}\fi  
   \ifx\thethanks@\empty\else \makefootnote@{}%
        {\thethanks@}\fi   
 \fi
  \pretitle
  \begingroup 
  \ifmonograph@ \topskip7pc \else \topskip4pc \fi
  \box\titlebox@
  \endgroup
  \preauthor
  \ifvoid\authorbox@\else \vskip2.5pcplus1pc\unvbox\authorbox@\fi
  \preaffil
  \ifvoid\affilbox@\else \vskip1pcplus.5pc\unvbox\affilbox@\fi
  \predate
  \ifx\thedate@\empty\else
       \vskip1pcplus.5pc\line{\hfil\thedate@\hfil}\fi
  \preabstract
  \ifvoid\abstractbox@\else
       \vskip1.5pcplus.5pc\unvbox\abstractbox@ \fi
  \ifvoid\tocbox@\else\vskip1.5pcplus.5pc\unvbox\tocbox@\fi
  \prepaper
  \vskip2pcplus1pc\relax}
\def\foliofont@{\tenrm} 
\def\headlinefont@{\tenpoint} 


\def\refstyle#1{\uppercase{%
  \if#1A\relax \def\keyformat##1{[##1]\enspace\hfil}%
  \else\if#1B\relax 
      \def\keyformat##1{\aftergroup\kern
           aftergroup-\aftergroup\refindentwd}%
      \refindentwd\parindent 
  \else\if#1C\relax 
      \def\keyformat##1{\hfil##1.\enspace}%
  \else\if#1D\relax 
      \def\keyformat##1{\hfil\llap{{##1}\enspace}} 
  \fi\fi\fi\fi}
}
\refstyle{D}   
\def\refsfont@{\tenpoint} 
\def\no{\gdef\no{\makerefbox\no\keybox@\empty}\no} \newbox\keybox@
\def\email{\let\savedef@\email
  \def\email##1\endemail{\let\email\savedef@
  \toks@{\def\usualspace{{\it\enspace}}\endgraf\indent\tenpoint}%
  \toks@@{\tt ##1\par}%
  \expandafter\xdef\csname email\number\addresscount@\endcsname
  {\the\toks@\frills@{{\noexpand\smc E-mail address\noexpand\/}:%
     \noexpand\enspace}\the\toks@@}}%
  \nofrillscheck\email}

\chardef\oldatsign=\catcode`\@
\catcode`\@=11
\newif\ifdraftmode 
\global\draftmodefalse


%
%
%
%
\newif\ifl@beloutopen
\newwrite\l@belout
\newread\l@belin

\global\let\currentfile=\jobname

\def\getfile#1{%
        \immediate\closeout\l@belout
        \global\l@beloutopenfalse
        \gdef\currentfile{#1}%
        \input #1%
        \par
        \newpage
}

\def\getxrefs#1{%
        \bgroup
                \def\gobble##1{}
                \edef\list@{#1,}%
                \def\gr@boff##1,##2\end{
                        \openin\l@belin=##1.ref
                        \ifeof\l@belin
                        \else
                                \closein\l@belin
                                \input ##1.ref
                        \fi
                        \def\list@{##2}%
                        \ifx\list@\empty
                                \let\next=\gobble
                        \else
                                \let\next=\gr@boff
                        \fi
                        \expandafter\next\list@\end
                }%
                \expandafter\gr@boff\list@\end
        \egroup
}

\def\testdefined#1#2#3{%
        \expandafter\ifx
        \csname #1\endcsname
        \relax
        #3%
        \else #2\fi
}

\def\document{%
        \minaw@11.11128\ex@ 
        \def\alloclist@{\empty}%
        \def\fontlist@{\empty}%
        \openin\l@belin=\jobname.ref 
        \ifeof\l@belin\else
                \closein\l@belin
                \input \jobname.ref
        \fi
}

\def\getst@te#1#2{%
        \edef\st@te{\csname #1s!#2\endcsname}%
        \expandafter\ifx\st@te\relax
                \def\st@te{0}%
        \fi
}

\def\setst@te#1#2#3{%
        \expandafter
        \gdef\csname #1s!#2\endcsname{#3}%
}

\outer\def\setupautolabel#1#2{%
        \def\newcount@{\global\alloc@0\count\countdef\insc@unt} 
        \def\newtoks@{\global\alloc@5\toks\toksdef\@cclvi}
        \expandafter\newcount@\csname #1Number\endcsname
        \expandafter\global\csname #1Number\endcsname=1%
        \expandafter\newtoks@\csname #1l@bel\endcsname
        \expandafter\global\csname #1l@bel\endcsname={#2}%
}

\def\reflabel#1#2{%
        \testdefined{#1l@bel}
        {
                \getst@te{#1}{#2}%
                \ifcase\st@te
                        ?
                        \message{Unresolved forward reference to
                                label #2. Use another pass.}%
                \or     
                        \setst@te{#1}{#2}2
                        \csname #1l!#2\endcsname 
                \or     
                        \csname #1l!#2\endcsname 
                \or     
                        \csname #1l!#2\endcsname 
                \fi
        }{
                {\escapechar=-1 
                \errmessage{You haven't done a
                        \string\\setupautolabel\space for type #1!}%
                }%
        }%
}

{\catcode`\{=12 \catcode`\}=12
        \catcode`\[=1 \catcode`\]=2
        \xdef\Lbrace[{] 
        \xdef\Rbrace[}]%
]%

\def\setlabel#1#2{%
        \testdefined{#1l@bel}
        {
                \edef\templ@bel@{\expandafter\the
                        \csname #1l@bel\endcsname}%
                \def\@rgtwo{#2}%
                \ifx\@rgtwo\empty
                \else
                        \ifl@beloutopen\else
                                \immediate\openout\l@belout=\currentfile.ref
                                \global\l@beloutopentrue
                        \fi
                        \getst@te{#1}{#2}%
                        \ifcase\st@te
                        \or     
                        \or     
                                \edef\oldnumber@{\csname #1l!#2\endcsname}%
                                \edef\newnumber@{\templ@bel@}%
                                \ifx\newnumber@\oldnumber@
                                \else
                                        \message{A forward reference to label
                                                #2 has been resolved
                                                incorrectly.  Use another
                                                pass.}%
                                \fi
                        \or     
                                \errmessage{Same label #2 used in two
                                        \string\setlabel s!}%
                        \fi
                        \expandafter\xdef\csname #1l!#2\endcsname
                                {\templ@bel@}
                        \setst@te{#1}{#2}3%
                        \immediate\write\l@belout 
                                {\string\expandafter\string\gdef
                                \string\csname\space #1l!#2%
                                \string\endcsname
                                \Lbrace\templ@bel@\Rbrace
                                }%
                        \immediate\write\l@belout 
                                {\string\expandafter\string\gdef
                                \string\csname\space #1s!#2%
                                \string\endcsname
                                \Lbrace 1\Rbrace
                                }%
                \fi
                \templ@bel@ 
                \expandafter\ifx\envir@end\endref 
                        \gdef\marginalhook@{\marginal{#2}}%
                \else
                        \marginal{#2}
                \fi
                \expandafter\global\expandafter\advance 
                        \csname #1Number\endcsname
                        by 1 %
        }{
                {\escapechar=-1
                \errmessage{You haven't done a \string\\setupautolabel\space
                        for type #1!}%
                }%
        }%
}


\newcount\SectionNumber
\setupautolabel{t}{\number\SectionNumber.\number\tNumber}
\setupautolabel{r}{\number\rNumber}
\setupautolabel{T}{\number\TNumber}

\define\rref{\reflabel{r}}

\define\rnum{\setlabel{r}}


\def\D{\text{\rm D}}

\def\F{{\Bbb F}}

\def\M{\text{\rm M}}

\def\PGL{\operatorname{PGL}}

\def\GL{\operatorname{GL}}
\def\Aut{\operatorname{Aut}}

\def\Nor{\text{\cal N}}

\def\proof{\noindent{\bf Proof.\ }}


\newcount\secno

\secno=1

\topmatter

\title
A note on the arrangement of subgroups in the\\
automorphism groups of submodule lattices of\\
free modules
\endtitle

\author
Alexandre A.$\,$Panin
\endauthor

\affil
{\it Department of Mathematics and Mechanics\\
St.Petersburg State University\\
2 Bibliotechnaya square,\\
St.Petersburg 198904, Russia}
\endaffil



\thanks
This research has been carried out in the framework of the
program ``Young Scientists''--1999.
\endthanks

\address
\endaddress
\email
alex@ap2707.spb.edu
\endemail

\date
\enddate

\abstract
A complete description of subgroups in the general linear group over a
semilocal ring containing the group of diagonal matrices was obtained
by Z.I.Borewicz and N.A.Vavilov. It is shown in the present paper that
a similar description holds for the intermediate subgroups of the
group of all automorphisms of the lattice of right submodules of a
free finite rank $R$--module over a simple Artinian ring containing
the group consisting of those automorphisms which leave invariant an
appropriate sublattice. Bibliography: 12 titles.
\endabstract

\endtopmatter

\document

\heading
\S~\the\secno. Introduction
\endheading
\advance\secno by 1

The following result is due to Z.I.Borewicz and N.A.Vavilov (see [\rref{BV}]
for a preliminary version).

\proclaim
{Theorem [\rref{V1}]} Let $R$ be a semilocal ring such that the decomposition
of the factor-ring $R/J(R)$ in the direct sum of simple Artinian rings
does not include either fields containing less than seven elements, or
the full matrix ring $\M(2,\F_2)$. Then for every intermediate
subgroup $F,\ \D(n,R)\leqslant F\leqslant \GL(n,R)$, there exists a
unique $D$--net $\sigma$ of two-sided ideals in $R$ of order $n$ such
that $G(\sigma)\leqslant F\leqslant \Nor(\sigma)$, where
$\Nor(\sigma)$ is the normaliser of the net subgroup $G(\sigma)$ in
$\GL(n,R)$.
\endproclaim

One may find a wealth of background information and many further
related references in the surveys [\rref{V2}], [\rref{V3}].

This theorem was generalised by A.Z.Simonian [\rref{S}], A.A.Panin and
A.V.Yakov\-lev~[\rref{PY}], A.A.Panin [\rref{P1}], [\rref{P2}]. It was shown
that these results can be stated in terms of a Galois correspondence between
sublattices of the lattice of (right) submodules of a free module
$R^n$ and subgroups of $\GL(n,R)$ considered as automorphism groups of
this lattice.

One can try to describe subgroups of the group of {\it all}
automorphisms of this lattice containing the group of ``diagonal''
automorphisms (i.e., the group consisting of those automorphisms which
leave invariant its appropriate sublattice).

We solve this problem for a simple Artinian ring in the present paper.

\heading
\S~\the\secno. Preliminaries
\endheading
\advance\secno by 1

Let $L$ be a lattice and let $G$ be a subgroup of the group $\Aut(L)$
of all automorphisms of the lattice $L$. Consider a subgroup $F$ of the
group $G$ and a sublattice $M$ of the lattice $L$. By definition, put

\medbreak

$L(F)=\{l\in L$ such that $f(l)=l$ for every $f\in F\}$,

\smallskip

$G(M)=\{g\in G$ such that $g(m)=m$ for every $m\in M\}$.

\medbreak

Let $L$ be a complete modular lattice and $L_0$ be a sublattice of $L$
which is a finite Boolean algebra such that $1_L=1_{L_0}$ and
$0_L=0_{L_0}$. Consider a subgroup $G$ of the group $\Aut(L)$. Let
$H=G(L_0),\ L'_0=L(H)$. By $e_1,e_2,\ldots,e_n$ we denote the
atoms of $L_0$.

For every $x\in L$ the collection $([x]_1,\ldots,[x]_n)$, where
$[x]_i=(x+\sum\limits_{j\neq i}e_j)\cdot e_i$, is called the {\it
support} of $x$ and is denoted by $[x]$. The support can be defined
equivalently as the minimal (with respect to the natural ordering)
collection $(x_1,\ldots,x_n)$, where $x_i\leqslant e_i$ and
$x\leqslant\sum\limits_{i=1}^nx_i$ (see [\rref{P2}], [\rref{PY}]).

Certain subgroups of type $G(M)$, where $M$ is a sublattice of $L'_0$,
are of special interest. The following definition was introduced in
[\rref{PY}].

\proclaim
{Definition}
A {\it net collection in} $L'_0$
is a collection of elements $\tau=(\tau_{ij})_{i,j=1}^n$ such
that for every indices $i,j,k$ and $g\in G$:
\smallskip
$(i)\ \tau_{ij}\leqslant e_j$
\smallskip
$(ii)\ \tau_{ii}=e_i$
\smallskip
$(iii)\ \tau_{ij}\in L'_0$
\smallskip
$(iv)\ [g(e_i)]_j\leqslant \tau_{ij}$ if and only if \
$[g(\tau_{ki})]_j\leqslant\tau_{kj}.$
\endproclaim

For every net collection $\tau=(\tau_{ij})$ in $L'_0$ the lattice
generated by elements $\sum\limits_{j=1}^{n}\tau_{ij},\ i=1,\ldots,n$
is called the lattice {\it associated} with $\tau$ and is denoted by
$K_{\tau}$.

\heading
\S~\the\secno. The main theorem
\endheading
\advance\secno by 1

Let $L$ be the lattice of right submodules of a free $R$--module $M$
of rank $n\geqslant 2$, where $R$ is a semilocal ring such that the field
$\F_2$ of two elements does not occur in the decomposition of the factor-ring
$R/J(R)$ in the direct sum of simple Artinian rings. Consider a basis
$\bar{e}_1,\ldots,\bar{e}_n$ of $M$ and let $e_i$ be spanned by
$\bar{e}_i$. Let $L_0$ be the sublattice of $L$ generated by
$e_1,\ldots,e_n$. One can consider $\PGL(n,R)$ as a subgroup of
$\Aut(L)$. Since we are only interested in description of the {\it
intermediate} subgroups, it does not make a difference to speak about
$\PGL(n,R)$ or $G=\GL(n,R)$ here, with $H=G(L_0)$ playing the role of
$\D(n,R)$ in our case.

It is clear that $\sigma$ is a net collection in
$L'_0$ if and only if the transposed collection of ideals
$\sigma^{\hbox{\small T}}$ is a $D$--net over $R$, and in this
case the net subgroup $G(\sigma^{\hbox{\small T}})$ equals $G(K_{\sigma})$.

\smallskip

The main result of this article is the following

\proclaim
{Theorem 1} Let $T$ be a skew field, $R=\M(m,T)$ a full matrix ring over
$T$ not equal to $\F_2,\F_3,\F_4,\F_5,\M(2,\F_2)$, and $G=\Aut(L)$ the group
of all automorphisms of $L$. Let $H=G(L_0)$, $L'_0=L(H)$. Then for
every subgroup $F,\ H\leqslant F\leqslant G$, there exists a unique
$D$--net $\sigma$ of two-sided ideals in $R$ of order $n$ such that
$G(K_{\sigma})\leqslant F\leqslant \Nor_GG(K_{\sigma})$, where
$K_{\sigma}$ is the lattice associated with the net collection
$\sigma^{\hbox{\small T}}$ in $L'_0$.
\endproclaim

\proof First, note that $\sigma^{\hbox{\small T}}$ is still
a net collection in $L'_0$ for every $D$--net $\sigma$ over $R$.

It is clear that the lattice of right $\M(m,T)$--submodules
of $\M(m,T)^n$ is isomorphic to the lattice of right $T$--subspaces
of $T^{nm}$. 

Let $n\geqslant 3$. Using the Fundamental Theorem of
Projective Geometry [\rref{A}], we can assume that
$G=\GL(n,R)\leftthreetimes\Aut(T)$
and $H=\D(n,R)\leftthreetimes\Aut(T)$. Thus each intermediate subgroup of
$G$ containing $H$ is of the form $F\leftthreetimes\Aut(T)$, where
$\D(n,R)\leqslant F\leqslant\GL(n,R)$. If the ring $R$
satisfies the conditions of Theorem~1, then for every such $F$ there
exists a unique $D$--net $\sigma$ of two-sided ideals in $R$ of order $n$ such
that $G(\sigma)\leqslant F\leqslant\Nor(\sigma)$. It is clear that
$G(K_{\sigma})=G(\sigma)\leftthreetimes\Aut(T)$. We have to prove
that $G(K_{\sigma})$ is normal in $F\leftthreetimes\Aut(T)$ (it is not
evident since $\Aut(T)$ is not a normal subgroup of $G$). If $f\in F$,
then $f^{-1}hf\in G(K_{\sigma})$ for every $h\in\D(n,R)$, therefore
$f(\sum\limits_{i=1}^n\sigma_{ij})\in L'_0$. Since $R$ is a simple ring,
we obtain $L'_0=L_0$. Hence $f^{-1}\varphi f\in G(K_{\sigma})$ for every
$\varphi\in\Aut(T)$. The rest is trivial.

The Fundamental Theorem of Projective Geometry fails for $n=2$, so one
has to consider this case separately. We can assume that $R=T$ is a skew field.
Then every proper nonzero subspace of $M$ is an atom of the lattice $L$,
so the automorphism group of $L$ is organised very simply: it's just the
symmetric group acting on the set of atoms of $L$. Let $T$ be a skew field
not equal to $\F_2$ and $\F_3$. Then $L$ contains at least 5 atoms. It's
easy to verify that there are only 5 intermediate subgroups: $H, \Nor_GH,
G(M_1),G(M_2), G$, where $M_1=\{e_1\},M_2=\{e_2\}$, and all subgroups except
$H$ are self-normalisable. This completes the proof of Theorem~1.

\medskip

If $n\geqslant 3$, then the assertion of Theorem~1 holds true for a somewhat
broader class of rings.

\smallskip

Let $R$ be a ring with the following property:

\smallskip

{\tt each automorphism of $R$ leaves invariant all its two-sided ideals}

\smallskip

Such objects will be called {\it good} rings. It is clear that full matrix
rings over skew fields, rings without automorphisms different from identical,
and uniserial rings are good ones (recall that a ring is called
uniserial, if its ideals are linearly ordered by inclusion).

\proclaim
{Theorem 2}
Let $R$ be a good semilocal ring such that the decomposition
of the factor-ring $R/J(R)$ in the direct sum of simple Artinian rings
does not include either fields containing less than seven elements, or
the full matrix ring $\M(2,\F_2)$. Let $n\geqslant 3$, $G=\Aut(L)$ be
the group of all automorphisms of $L$, $H=G(L_0)$, $L'_0=L(H)$.
Then for every subgroup $F,\ H\leqslant F\leqslant G$, there exists a unique
$D$--net $\sigma$ of two-sided ideals in $R$ of order $n$ such that
$G(K_{\sigma})\leqslant F\leqslant \Nor_GG(K_{\sigma})$, where
$K_{\sigma}$ is the lattice associated with the net collection
$\sigma^{\hbox{\small T}}$ in $L'_0$.
\endproclaim

\proof Follow the lines of the proof of Theorem~1, using the Fundamental
Theorem of Projective Geometry for rings of stable rank $2$ [\rref{Ve}].

\heading
\S~\the\secno. Final remarks
\endheading
\advance\secno by 1

$1^o.$~The Fundamental Theorem of Projective
Geometry holds in more general settings [\rref{OjS}], [\rref{SV}], [\rref{Ve}].
Therefore it is natural to try to generalise Theorems~1 and~2 in order to
cover semilocal rings (i.e. to formulate the result in terms of lattices and
their automorphism groups).

\smallskip

$2^o.$~The assertion of Theorem~1 for $n=2$ and 
skew fields can be deduced from a more general theorem on the
automorphisms of certain complete modular lattices, see [\rref{P2}].

\proclaim
{Definition}
For every $i\neq j$ and $x\leqslant e_j$
the set of $f\in G$ such that
\smallskip
$(i)\ f(e_s)=e_s$ for $s\ne i$,
\smallskip
$(ii)\ [f(e_i)]_k=\cases 0,&k\neq i,j;\\ e_i,&k=i;\\ x,&k=j\\ \endcases$
\smallskip\noindent
is denoted by $H_{ij}(x)$ and its elements are called ``transvections''.
\endproclaim

The following result is a corollary of Theorem~3.1~[\rref{P2}] (see
also [\rref{PY}] for a slightly weaker result).

\proclaim
{Theorem 3}
Let $L$ be a modular lattice of finite length, $L_0$ a
sublattice of the same length which is a Boolean algebra with atoms
$e_1,\ldots,e_n$, let $G$ be a subgroup of the group of all
automorphisms of the lattice $L,\ H=G(L_0)$.
\smallskip\indent
Assuming that the conditions $(a)-(d)$ stated below are satisfied, for
every subgroup $F\geqslant H$ of the group $G$ the group
$G(K_{\sigma})$ is normal in $F$.
\smallskip\indent
$(a)$~For every $i$ there exists at least one automorphism from $H$
which changes all atoms $x\in L \setminus \{ e_i \} $ such that
$[x]_i=e_i$ and leaves invariant all atoms $x\in L$ with $[x]_i=0$.
\smallskip\indent
$(b)$~If $x, y\in L$ are atoms with $[x]=[y]$, then there exists
$h\in H$ such that $h(x)=y$.
\smallskip\indent
$(c)$~For every $i\neq j$ the set $H_{ij}(e_j)$ is not
empty.
\smallskip\indent
$(d)$~For every $a\in G$ and $i\ne j$ the subgroup generated by $h\in
H$ such that $H_{ij}([aha^{-1}(e_i)]_j)\cap \langle
a,H\rangle\neq\varnothing$, and $H_{ij}([ah^{-1}a^{-1}(e_i)]_j)\cap
\langle a,H\rangle\neq\varnothing$ is equal to $H$.
\endproclaim

It is worth mentioning that the description of the intermediate
subgroups of the general linear group over a skew field, containing the
group of diagonal matrices, can also be deduced from this theorem. Moreover,
it is likely that methods and results of [\rref{P1}], [\rref{P2}], [\rref{PY}]
can bring further insight in our understanding of linear groups.

\heading
Acknowledgements
\endheading

I would like to thank Nikolai A.$\,$Vavilov and Elizaveta
V.$\,$Dybkova for their interest and helpful suggestions.

\enddocument